\documentclass[review]{elsarticle}

\usepackage{lineno,hyperref}
\modulolinenumbers[5]

\usepackage{stmaryrd}
\usepackage{amssymb}
\usepackage{color, amsmath,amssymb, amsfonts, amstext,amsthm, latexsym}

\usepackage{amssymb, epsfig, amssymb, latexsym}
\usepackage{amsmath}
\usepackage{graphicx}
\usepackage{longtable}
\usepackage{graphicx}
\usepackage{subfigure}
\usepackage{epstopdf}

\renewcommand{\phi}{\varphi}

\newtheorem{example}{Example}

\begin{document}

\begin{frontmatter}

\title{Variational inference of the drift function for stochastic differential equations driven by L\'{e}vy processes\tnoteref{mytitlenote}}


\author[mymainaddress]{Min Dai}
\ead{mindai@hust.edu.cn}

\author[myfourthaddress]{Jinqiao Duan}
\cortext[mycorrespondingauthor]{Corresponding author}
\ead{duan@iit.edu}

\author[mymainaddress]{Jianyu Hu}

\author[mymainaddress]{Xiangjun Wang}


\address[mymainaddress]{School of Mathematics and Statistics, \& Center for Mathematical Science, \\ Huazhong University of Science and Technology, Wuhan, 430074, China.}
\address[myfourthaddress]{Department of Applied Mathematics, College of Computing,\\ Illinois Institute of Technology, Chicago, IL 60616, USA.}

\begin{abstract}
  In this paper, we consider the nonparametric estimation problem of the drift function of stochastic differential equations driven by $\alpha$-stable L\'{e}vy motion. First, the Kullback-Leibler divergence between the path probabilities of two stochastic differential equations with different drift functions is optimized. By using the Lagrangian multiplier, the variational formula based on the stationary Fokker-Planck equation is constructed. Then combined with the data information, the empirical distribution is used to replace the stationary density, and the drift function is estimated non-parametrically from the perspective of the process. In the numerical experiment, the different amounts of data and different $\alpha$ values are studied. The experimental results show that the estimation result of the drift function is related to both. When the amount of data increases, the estimation result will be better, and when the $\alpha$ value increases, the estimation result is also better.
\end{abstract}


\end{frontmatter}


\section{Introduction}
Stochastic dynamical systems are often used to describe random phenomena in real life. Usually, stochastic differential equation\cite{Oksendal,Gihman} is used to model the stochastic process, but one of the important issues is the fitting of the model and the observed data, that is, the system identification in the data. The stochastic differential equation consists of a drift function and a diffusion part. For the thermodynamic equilibrium model, the diffusion function is proportional to the identity matrix, and the drift function is the gradient of potential energy. A simple and known method for estimating the drift function from the data is obtained from\cite{Iacus}, which takes advantage of the fact that the potential function is proportional to the logarithm of the stationary density. Therefore, we can get an explicit representation of the drift function by using the kernel density (KDE) estimator. However, this estimator is only based on the distribution of the data, completely ignoring the time sequence of the observations and the time lag between the data. Kernel density estimator is non-parameterized, it does not require the drift function satisfies the specific parameterized form. For non-equilibrium model, in general, the explicit representation between the drift function and the density is unknown. At the same time, for the higher dimensional case, the convergence of the kernel density to the real density will be slow with the increase of sample data\cite{Sriperumbudur}.

There are many methods for estimating the parameters of the drift function\cite{Batz,Ruttor}. For these parametric and non-parametric methods, we have to deal with the problem of low data sampling rate\cite{Gottschall}. For example, the non-parametric method of conditional expectation is used\cite{Honisch}, which requires numerical solution of the Kolmogorov backward equation at a given time interval, and also requires a large amount of data to calculate the conditional expectation. The Bayesian estimator using the Gaussian process prior gives a good estimate of the complete path under dense observations\cite{Papaspiliopoulos}. However, generally speaking, this method requires that unobserved diffusion paths be classified as hidden random variables between adjacent observations, which may lead to time-consuming calculations or require further approximation. Papaspiliopoulos et al.\cite{Papaspiliopoulos} introduced the Monte Carlo Gibbs sampling technique, which switched back and forth between hidden path sampling and drift function sampling of the process. For the processing of hidden paths, the expectation maximization method\cite{Rut} can be used to approximate the hidden process with linear stochastic differential equations. However, these methods only consider the stochastic differential equations driven by Brownian motion. In nature, there exists a class of stochastic differential equations driven by non-Gaussian L\'{e}vy motion. The L\'{e}vy motion\cite{Sato,Applebaum}, like the Brownian motion, is a type of Markov process with independent and stationary increments, which has a wide applications in the fields of physics, biology, earth science, etc.  Therefore, it is of great significance to study the data recognition problem of this kind of the stochastic differential equation. For the identification problem of stochastic differential equations driven by L\'{e}vy motion, Li and Duan\cite{Li} studied this equation by using parameterized identification method based on Kramers-Moyal coefficient. However, due to the need to solve the conditional expectation value, it requires a high amount of data and increases the calculation cost.

In this paper, the main purpose is to construct an effective nonparametric inference of the drift function of a class of stochastic differential equations driven by $\alpha$-stable L\'{e}vy motion, which can reduce the requirement for data sample size to a certain extent. In our method, we construct the Kullback-Leibler divergence\cite{Kullback} of the transition density through the diffusion process with two different drift functions, and use the variational formula based on the stationary Fokker-Planck equation\cite{Risken} to derive the explicit form of the drift function. Combined with data information, we use empirical distribution to replace stationary density, and further minimize empirical functional to obtain target results. However, if the empirical functional has a penalty function term based on the kernel, then the method can be extended to nonparametric inference. In addition, the properties of $\alpha$-stable L\'{e}vy motion will result in a large variation due to the different $\alpha$ values, we also carry out numerical experiments for different parameters $\alpha$.

This paper is arranged as follows. In section 2, We derive the variational formula of the Fokker-Planck equation corresponding to the stochastic differential equation driven by an $\alpha$-stable L\'{e}vy motion, and we use empirical distribution to replace stationary density and make full use of data information to give the expression of the drift function. In section 3, We show the feasibility of this variational method through numerical application examples. We end the paper some conclusions and discussions in section 4.

\section{Theory}
In this work, we consider the stochastic differential equation for the dynamics of a $d$-dimensional diffusion processes $X_t\in\mathbb{R}^d$ given by
\begin{equation}\label{eq:1}
dX_t=g(X_t)dt+\sigma(X_t){dB_t}+{dL_t^\alpha},
\end{equation}
where $g(\cdot)\in\mathbb{R}^d$ is the drift function, the diffusion function $\sigma(\cdot)$ is the $d\times{k}$ dimensional matrix, $B_t$ is the Brownian motion in $\mathbb{R}^k$ and $L_t^\alpha$ is a symmetric $\alpha$-stable L$\mathrm{\acute{e}}$vy motion in $\mathbb{R}^d$ with the generating triplet $(0,0,\nu_\alpha)$. The jump measure\cite{Duan}
\begin{equation*}
\nu_\alpha(dy)=c(n,\alpha)||y||^{-(n+\alpha)}dy,
\end{equation*}
with $c(n,\alpha)=\frac{\alpha\Gamma((n+\alpha)/2)}{2^{1-\alpha}\pi^{n/2}\Gamma(1-\alpha/2)}.$ The processes $B_t$ and $L_t^\alpha$ are taken to be independent.

\subsection{\textbf{Variational formulation for the Fokker-Planck equation}}
Now, we will discuss how can we determine the drift function $g$ from data of sample path under the noise matrix $\sigma(\cdot)$ known. Suppose that the stationary density $p(x)$ of the process are given. We then splits the drift into two parts $g(x)=r(x)+f(x)$, where $r(x)$ is known part and $f(x)$ is the part that we try to compute. Of course, in the multivariate case there is not enough information to reconstruct $f$ uniquely. However, we may search for a minimal solution which minimizes a quadratic functional
\begin{equation}\label{eq:2}
\frac{1}{2}\int{p(x)f(x)\cdot{A^{-1}(x)f(x)}}dx
\end{equation}
for a given positive definite matrix $A(x)$. Introducing a Lagrangian multiplier function $\psi(x)$ for the condition that the density $p$ fulfills the stationary Fokker-Planck equation with drift $g$, we can derive the minimal $f$ from the variation of the Lagrange-functional
\begin{equation}\label{eq:3}
\frac{1}{2}\int{p(x)f(x)\cdot{A^{-1}(x)f(x)}}dx-\int{\psi(x)\{\mathcal{L}p(x)-\nabla\cdot(f(x)p(x))\}}dx
\end{equation}
where the Fokker-Planck operator $\mathcal{L}$ corresponding to known drift $r(x)$ is given by
\begin{equation}\label{eq:4}
\mathcal{L}p(x)=-\nabla\cdot(r(x)p(x))+\frac{1}{2}{\mathrm{tr}[\nabla\nabla^T(D(x)p(x))]}+\int_{\mathbb{R}^d\setminus\{0\}}{[p(x+y,t)-p(x)]\nu_\alpha(dy)}
\end{equation}
with $D(x)\doteq{\sigma(x)\sigma(x)^T}.$ Variation of \eqref{eq:3} with respect to $f$ yields $f(x)=A(x)\nabla\psi(x).$ Inserting this solution back into \eqref{eq:3} shows that the unknown function $\psi$ can be derived from the minimization of the functional
\begin{equation}\label{eq:5}
\varepsilon[\psi]=\int\Bigg[\frac{1}{2}\nabla{\psi(x)}\cdot{A(x)\nabla{\psi(x)}}+\mathcal{L}^*\psi(x)\Bigg]p(x)dx,
\end{equation}
where $\mathcal{L}^*$ is the adjoint operator of $\mathcal{L}$, \eqref{eq:4} which fulfills
\begin{equation}\label{eq:6}
\int{\psi(x)\mathcal{L}p(x)dx}=\int{p(x)\mathcal{L}^*\psi(x)dx}
\end{equation}
and is given by
\begin{equation}\label{eq:7}
\mathcal{L}^*\psi(x)=r(x)\cdot\nabla(x)+\frac{1}{2}{\mathrm{tr}[(D(x)\nabla\nabla^T\psi(x))]}+\int_{\mathbb{R}^d\setminus\{0\}}{[\psi(x+y,t)-\psi(x)]\nu_\alpha(dy)}
\end{equation}
In fact, a direct minimization of \eqref{eq:5} with respect to $\psi$ yields
\begin{equation}\label{eq:8}
\mathcal{L}[\psi]p(x)=\mathcal{L}p(x)-\nabla\cdot(A(x)\nabla\psi(x)p(x))=0,
\end{equation}
which is the stationary Fokker-Planck equation corresponding to the density $p(x)$ and the drift $g(x)=r(x)+A(x)\nabla\psi(x).$ For the special case $D=A=I$ and $r=0$ the functional \eqref{eq:5} was introduced in the field of machine learning as a score-function for estimating $\mathrm{ln}{p(x)}$ up to a normalization constant\cite{Hyv}.

\subsection{\textbf{Minimizing the empirical functional}}
Our goal is to estimate $\psi$ from data by replacing the average over the stationary density $p(x)$ in the functional \eqref{eq:5} by empirical distribution
\begin{equation}\label{eq:9}
\hat{p}(x)=\frac{1}{n}\sum_{i=1}^{n}\delta(x-x_i)
\end{equation}
where $x_1,\cdots,x_n$ is a random, ergodic sample drawn from this density. An obvious possibility to construct estimators is to work with a parametric representation
\begin{equation}\label{eq:10}
\psi_\omega(x)=\sum_{k=1}^{K}\omega_k\phi_k(x)
\end{equation}
where the $\phi_k$ are a set of given basis functions. Then we can obtain the empirical constraint optimization problem which is minimize the following
\begin{equation}\label{eq:11}
\varepsilon_{\mathrm{emp}}[\psi_\omega]=\sum_{i=1}^{n}\Bigg[\frac{1}{2}\nabla{\psi_{\omega}(x_i)}\cdot{A(x_i)\nabla{\psi_{\omega}(x_i)}}+\mathcal{L}^*\psi_{\omega}(x_i)\Bigg]
\end{equation}
which is a quadratic form in the $\omega_k$ and can thus be performed in closed form. We are however interested in the case where a representation in terms of a finite set of basis functions is not rich enough to represent $\psi$. Thus we will resort to a more general, nonparametric representation allowing for an infinite set of functions $\phi_k.$ So, we can understand \eqref{eq:10} as a Gaussian process model\cite{Rasmussen} for the function $\psi,$ that is we will choose the quadratic form $\frac{1}{2}\Sigma_k{\omega_k^2/\lambda_k}$ as an extra penalty term, where the $\lambda_k$ are hyper-parameters. The penalty can also be viewed from a pseudo-Bayesian perspective where $\mathrm{exp}{\{-C\varepsilon_{\mathrm{emp}}{[\psi_\omega]}\}}$ is interpreted as a likelihood and $\mathrm{exp}{\{-\frac{1}{2}\Sigma_k{\omega_k^2/\lambda_k}\}}$ as a Gaussian prior distribution over parameters $\omega_k$. $C$ can be chosen to give different weight to the data and to the penalty.

Motivated by Gaussian process point of view we will introduce the kernel trick into our formalism avoiding an explicit specification of $\phi_k$ and $\lambda_k$ and assume instead that these are defined implicitly as orthonormal eigenfunctions and eigenvalues of a positive definite kernel function $K(x,x^{'})$ via
\begin{equation*}
K(x,x^{'})=\sum_k{\lambda_k\phi_k(x)\phi_k(x^{'})}.
\end{equation*}
In the kernel approach the regularized functional can be written as
\begin{equation}\label{eq:12}
C\sum_{i=1}^{n}\Bigg[\frac{1}{2}\nabla{\psi(x_i)}\cdot{A(x_i)\nabla{\psi(x_i)}}+\mathcal{L}^*\psi(x_i)\Bigg]+\frac{1}{2}\int\int\psi(x)K^{-1}(x,x^{'})\psi(x^{'})dxdx^{'},
\end{equation}
where $K^{-1}(x,x^{'})$ is the formal inverse of the kernel operator. Performing the variation with respect to $\psi$ yields
\begin{equation*}
Cn\mathcal{L}[\psi]\hat{p}(x)+\int{K^{-1}(x,x^{'})\psi(x^{'})dx^{'}}=0,
\end{equation*}
where $\mathcal{L}[\psi]$ was defined in \eqref{eq:8}. Multiplying both sides of this equation with the operator $K$ we get
\begin{equation}\label{eq:13}
\psi(x)+C\sum_{j=1}^{n}{\mathcal{L}_{x^{'}}^{*}[\psi]K(x,x^{'})_{x^{'}=x_j}}=0
\end{equation}
where the adjoint operator acts on functions $h$ as
\begin{equation}\label{eq:14}
{\mathcal{L}_{x^{'}}^{*}h(x^{'})}=(r(x^{'})+A(x^{'})\nabla\psi(x^{'}))\nabla h(x^{'})+\frac{1}{2}\mathrm{tr}{[D(x^{'})\nabla\nabla^Th(x^{'})]}+\int_{\mathbb{R}^d\setminus\{0\}}{[h(x^{'}+y)-h(x^{'})]\nu_\alpha(dy)}
\end{equation}
Applied the family of kernel functions $h(x^{'})=K(x,x^{'})$ when the stationary density $p$ is replaced by its empirical approximation $\hat{p}$ \eqref{eq:9}.

The gradient of $\psi$ at the data points is computed by taking the gradient of \eqref{eq:13} and setting $x=x_i.$ This yields the set of linear equations
\begin{equation}\label{eq:15}
\nabla\psi(x_i)+C\sum_{j=1}^{n}{\mathcal{L}_{x^{'}}^{*}[\psi]\nabla_xK(x,x^{'})_{x=x_i,x^{'}=x_j}}=0
\end{equation}
for the $d\times{n}$ unknowns $\nabla\psi(x_i)$ which can be plugged into \eqref{eq:13} to obtain the explicit result for the estimator.

\section{Application}
We will give a simple numerical example in this section to illustrate the results through the use of different amounts of data and the comparison of different $\alpha$ values.

\begin{example}
Consider a stochastic dynamical system of the following form
\begin{equation}\label{eq:16}
dX_{t}=(X_t-X_t^3)dt+{dW_{t}}+{dL_{t}^\alpha},~~~  X_0=x\in{\mathbb{R}},
\end{equation}
Where $W_t$ is the standard Brownian motion, $L_t^\alpha$ is a stable L\'{e}vy motion, the drift function is $f(x)=x-x^3$. In the numerical simulation, we obtain a sample path data of the time series by observing the process $X_t$ once. In this example, we take the time step $\Delta{t}=0.001$, and the sample sizes are $n=6000$ and $n=10000$ data points respectively, and then use the variational method of the Fokker-Planck equation to perform nonparametric inferences on the drift function $f$. In addition, we also select different $\alpha$ values to compare the results. we will show the numerical results:
\begin{figure}
\begin{minipage}[]{0.5 \textwidth}
\centerline{\includegraphics[width=8cm,height=6cm]{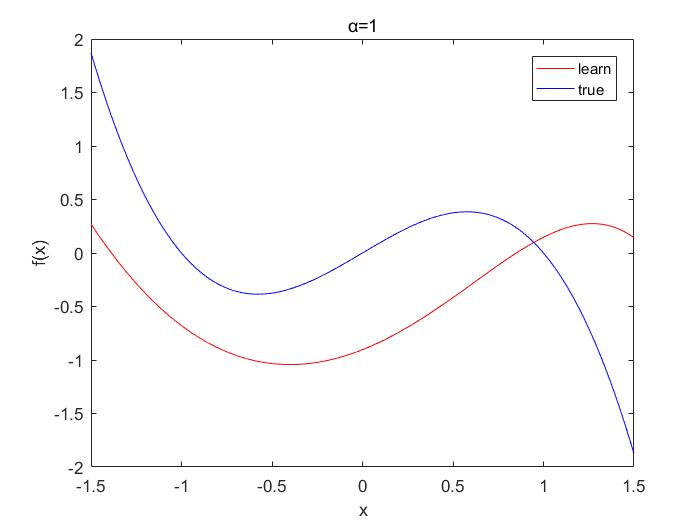}}
\end{minipage}
\hfill
\begin{minipage}[]{0.5 \textwidth}
\centerline{\includegraphics[width=8cm,height=6cm]{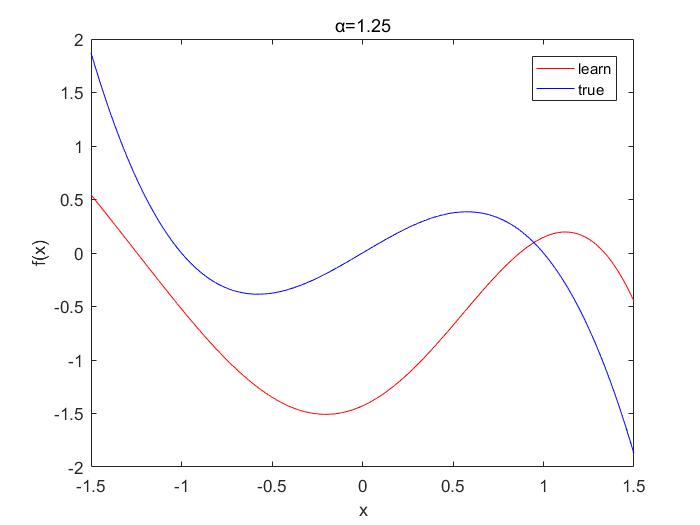}}
\end{minipage}
\hfill
\begin{minipage}[]{0.5 \textwidth}
\centerline{\includegraphics[width=8cm,height=6cm]{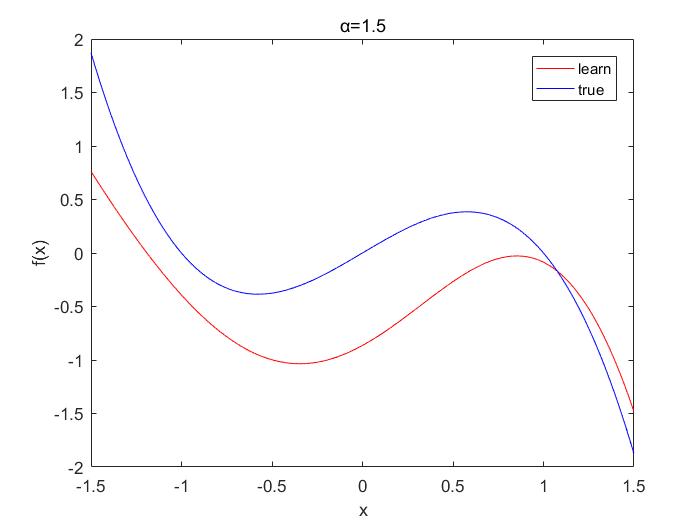}}
\end{minipage}
\hfill
\begin{minipage}[]{0.5 \textwidth}
\centerline{\includegraphics[width=8cm,height=6cm]{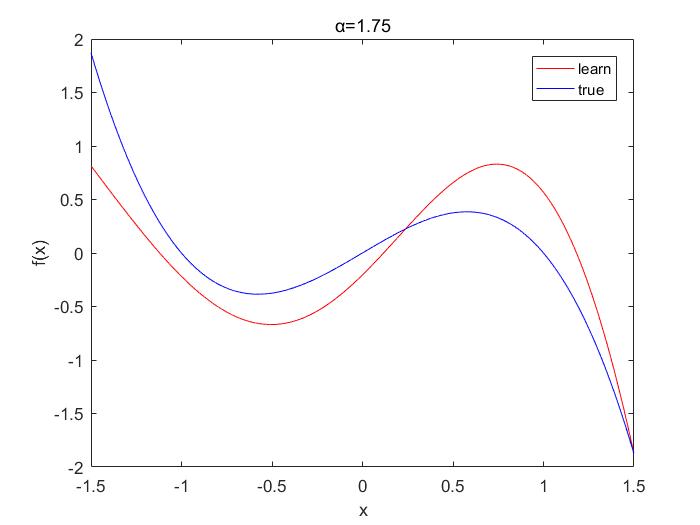}}
\end{minipage}
\caption{When the sample data size $n=6000$, the nonparametric inference results of the drift function $f$ with different $\alpha$ values.}
\label{f1-levy}
\end{figure}

\begin{figure}
\begin{minipage}[]{0.5 \textwidth}
\centerline{\includegraphics[width=8cm,height=6cm]{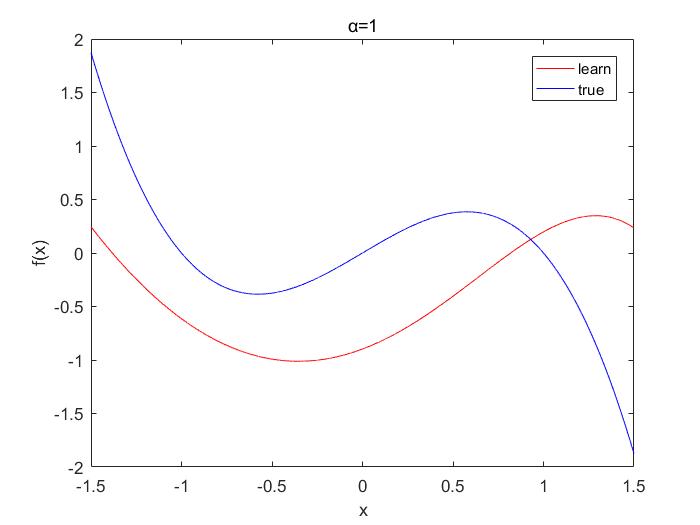}}
\end{minipage}
\hfill
\begin{minipage}[]{0.5 \textwidth}
\centerline{\includegraphics[width=8cm,height=6cm]{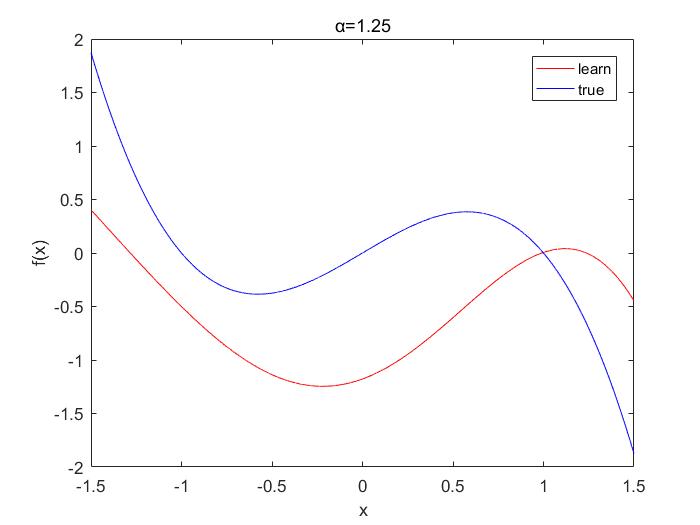}}
\end{minipage}
\hfill
\begin{minipage}[]{0.5 \textwidth}
\centerline{\includegraphics[width=8cm,height=6cm]{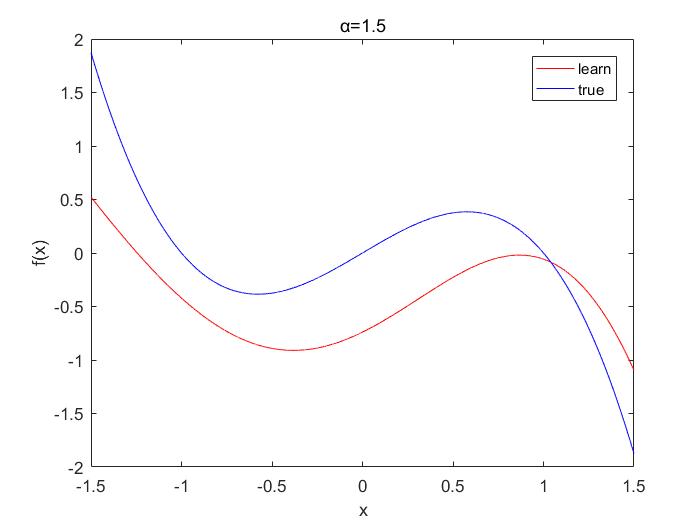}}
\end{minipage}
\hfill
\begin{minipage}[]{0.5 \textwidth}
\centerline{\includegraphics[width=8cm,height=6cm]{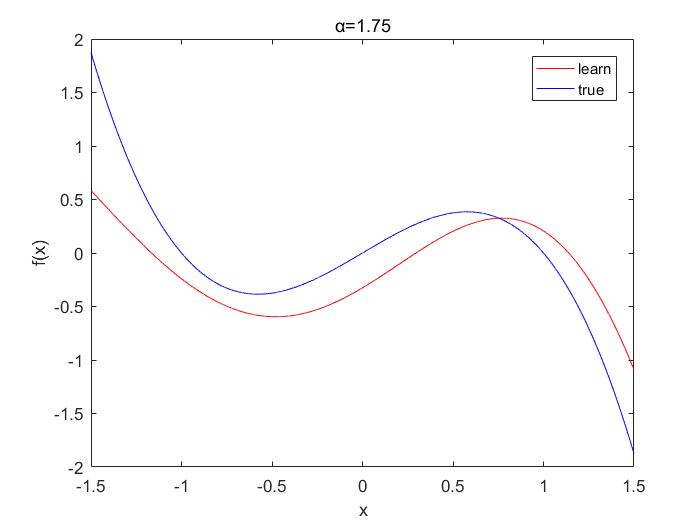}}
\end{minipage}
\caption{When the sample data size $n=10000$, the nonparametric inference results of the drift function $f$ with different $\alpha$ values.}
\label{f2-levy}
\end{figure}
\end{example}
From Figure \eqref{f1-levy} and \eqref{f2-levy}, we can see that the size of the data and the different $\alpha$ values ??have an impact on the inference of the drift function. In this chapter, our method only considers the data obtained from one observation of the process. When the value of $\alpha$ is the same, the estimator of the drift function becomes better as the amount of data increases.
And for the same amount of data, we're going to consider the selection of $\alpha$ values of 1, 1.25, 1.5 and 1.75, respectively. As the value of $\alpha$ increases, the sample trajectory of the process is smoother, making the inference result of the drift function better. Of course, we also know from the numerical results that the estimation of the drift function needs to be improved, such as when $\alpha$ is close to 1, the estimation effect of the drift function is poor, so we can improve the numerical results by increasing the amount of sample data. In addition, the method to infer the drift function only. For diffusion function, we are still in the further study. On the one hand, we hope to obtain the diffusion function inference in the case of additive and multiplicative noise when the drift function is known. On the other hand, we also hope to obtain an optimization method to infer the drift function and the diffusion function at the same time.

\section{Conclusion and discussions}
In this paper, we have proposed a way to obtain the drift function estimator for a class of stochastic differential equations driven by $\alpha$-stable L\'{e}vy motion from data. We construct a variational formula based on the stationary Fokker-Planck equation, replace the stationary density with an empirical distribution, and derive the drift function expression. Another contribution is that our method only needs to observe the dense data of the process once, which reduces the cost of data processing, and is more in line with the actual observation situation, so it is more conducive to processing the actual data.

In addition, there are some interesting extensions for this study. We only studied the nonparametric estimation of the drift function of the stochastic dynamical system driven by L\'{e}vy motion. It is worth thinking about the inference of the diffusion function for additive and multiplicative cases. On the other hand, we can also explore the method of extracting the dynamics of this type of the stochastic system from the data, so as to gain insight into the stochastic phenomenon in related fields.

\section*{Acknowledgements}
We would like to thank Dr. Xiaoli Chen, Dr. Wei Wei, Dr. Cheng Fang for helpful discussions. This work was supported by the National Natural Science Foundation of China (NSFC) grants 12001213, 11801192, 11771449 and National Science Foundation (NSF) grant 1620449.

\section*{Appendix: Kullback-Leibler divergence and Expression of estimators}
\medskip
\noindent \textbf{A1.  Kullback-Leibler divergence}

For every $t\in[0,T]$, assuming that the dense path observation $X_{0:T}$ of equation \eqref{eq:1} is obtained, we can derive the likelihood function. We discretize the time interval $[0,T]$, take the time step as $\Delta{t}$ and use the fact that $\Delta{t}\rightarrow{0}$, then the transition density of the equation \eqref{eq:1} satisfies the Gaussian distribution. We find that the negative log-likelihood part that depends on the drift function $g$ can be approximated by the following equation
\begin{equation*}
-\ln p(X_{0:T}|g)\simeq{\frac{1}{2}\sum_{t}\{\|g(X_t)\|^2\Delta{t}-2\langle{g(X_t),(X_{T+\Delta{t}}-X_t)}\rangle\}}+const,
\end{equation*}
Where the inner product $\langle{u,v}\rangle\doteq{u\cdot{D^{-1}}v}$, and the corresponding square norm is $\|u\|^2\doteq{u\cdot{D^{-1}}u}$.

For a sufficiently small $\Delta{t}$, using the Gaussian form of the transition density, we can give the Kullback-Leibler divergence between the path probabilities of the two diffusion processes with drift functions of $g(x)$ and $r(x)$ respectively, where $g(x)=f(x)+r(x)$, which can be expressed as
\begin{equation}\label{eq:17}
D(p(X_{0:T}|g)\|p(X_{0:T}|r))=\int_0^Tdt\int{p_t(x)f(x)\cdot{D^{-1}(x)f(x)dx}}.
\end{equation}
Assume that the two diffusion processes have the same diffusion function $D(x)$ and non-stochastic initial state \cite{Archambeau}. The density $p_t(x)$ is the probability density of the diffusion process with drift function $g$ at time $t$. Furthermore, we assume that the process reaches a stationary state and has a stationary density $p(x)$, thus the relative entropy rate \cite{Chernyak} can be obtained
\begin{equation}\label{eq:18}
\lim_{T\rightarrow{\infty}}\frac{1}{T}D(p(X_{0:T}|g)\|p(X_{0:T}|r))=\int{p_t(x)f(x)\cdot{D^{-1}(x)f(x)dx}}.
\end{equation}
The comparison with equation \eqref{eq:3} shows that when $A=D$, the minimization of equation \eqref{eq:5} will result in that the stationary density process with drift function $g$ is close enough to the process with drift function $r$ in the case of relative entropy. Therefore, we can understand this as the generalized minimum relative entropy (Kullback Leibler divergence) under the constraints given by the stationary density.

\medskip
\noindent \textbf{A2.  Expression of estimators}

We will show how to obtain an explicit representation of the estimator from equations \eqref{eq:13} and \eqref{eq:15}. In order to simplify the notation, we make some notation simplification for the expressions covered in the equation. Let the vectors of $\nabla\psi(x_i)$ and drift function $r(x_i)$ at all data points be
\begin{equation}\label{eq:19}
\begin{split}
&b=(\nabla\psi(x_1),\nabla\psi(x_2),\cdots,\nabla\psi(x_n)),\\
&r=(r(x_1),r(x_2),\cdots,r(x_n)).
\end{split}
\end{equation}
The kernel function $B$ and the block-diagonal weight matrix $A$ are defined as
\begin{equation}\label{eq:20}
\begin{split}
&B_{ij}=\nabla_{x}\nabla_{x'}K(x,x')|_{x=x_i, x'=x_j},\\
&A_{ij}=\delta_{ij}A(x_j),
\end{split}
\end{equation}
and define
\begin{equation}\label{eq:21}
\begin{split}
&y(x)=\frac{1}{2}\sum_{j=1}^ntr[D(x')\nabla_{x'}\nabla_{x'}^TK(x,x')]_{x'=x_j},\\
&z(x)=\int_{\mathbb{R}^d-\{0\}}[K(x,x'+w)-K(x,x')]\nu_{\alpha}(dw),
\end{split}
\end{equation}
then we have
\begin{equation}\label{eq:22}
\begin{split}
&\nabla y(x)=(\nabla{y(x_1)},\nabla{y(x_2)},\cdots,\nabla{y(x_n)}),\\
&\nabla z(x)=(\nabla{z(x_1)},\nabla{z(x_2)},\cdots,\nabla{z(x_n)}).
\end{split}
\end{equation}
We can rewrite \eqref{eq:15} as
\begin{equation}\label{eq:23}
(I+CBA)b+C(Br+\nabla y+\nabla z)=0,
\end{equation}
Where $I$ represents the identity matrix. Equation \eqref{eq:13} can be similarly written as
\begin{equation}\label{eq:24}
\psi(x)+C(k(x)Ab+k(x)r+y(x)+z(x))=0,
\end{equation}
and $k(x)=(\nabla_{x'}K(x,x')|_{x'=x_1},\nabla_{x'}K(x,x')|_{x'=x_2},\cdots,\nabla_{x'}K(x,x')|_{x'=x_n}).$ Further, we can obtain the following result via solving the function $b$ of equation \eqref{eq:23},
$$b=-(I+CBA)^{-1}C(Br+\nabla y+\nabla z).$$
Combined with equation \eqref{eq:24}, the explicit form of the estimator $\psi$ can be solved
\begin{equation}
\psi(x)=Ck(x)A(I+CBA)^{-1}C(Br+\nabla y+\nabla z)-C(k(x)r+y(x)+z(x)).
\end{equation}

\end{document}